\documentclass[11pt]{amsart}

\usepackage{mathptmx}
\usepackage{eucal}
\usepackage{amsmath}
\usepackage{amscd}
\usepackage{amssymb}
\usepackage{amsthm}
\usepackage{xspace}
\usepackage[all,tips]{xy}
\usepackage[dvips]{graphicx}
\usepackage{verbatim}
\usepackage{syntonly}
\usepackage{hyperref}

\theoremstyle{plain}
\newtheorem{thm}{Theorem}[section]

\newtheorem{lemma}[thm]{Lemma}

\theoremstyle{definition}
\newtheorem*{rem}{Remark}

\DeclareMathOperator{\CR}{cr}

\newcommand{\set}[1]{\left\{#1\right\}}

\newcommand{\pr}[1]{\left( #1 \right) }

\def\free{\mathbb{F}}
\def\into{\hookrightarrow}
\def\onto{\twoheadrightarrow}
\def\co{\colon\thinspace}
\newcommand{\mnote}[1]{}

\usepackage{fancyhdr}

\fancypagestyle{plain}

\title{\textbf{Graphs of subgroups of free groups}}
\author{Larsen Louder and D. B. McReynolds}
\thanks{Both authors were supported in part by NSF postdoctoral fellowships.}

\begin{document}

%---------------------------------------------------------------------------
%---------------------------------------------------------------------------

\begin{abstract}
  \noindent We construct an efficient model for graphs of finitely
  generated subgroups of free groups. Using this we give a very short
  proof of Dicks's reformulation of the strengthened Hanna Neumann
  Conjecture as the Amalgamated Graph Conjecture. In addition, we
  answer a question of Culler and Shalen on ranks of intersections in
  free groups. The latter has also been done independently by
  R. P. Kent IV.
\end{abstract}

\maketitle

%---------------------------------------------------------------------------
%---------------------------------------------------------------------------
\section{Introduction}

\thispagestyle{empty}

\noindent One purpose of this article is to investigate the interplay
between the join and intersection of a pair of finitely generated
subgroups of a free group. Our main result,
Theorem~\ref{reducetovalencethree}, is a minor generalization of the
construction of the first author from~\cite{louder}, and produces a
simple model for analyzing intersections and joins. We use this
technique to give a quick proof of a theorem of Dicks
\cite{dicks}. Another application of
Theorem~\ref{reducetovalencethree} is an answer to an
unpublished question of Culler and Shalen \cite{cullershalen}. This
has been done independently by Kent \cite{kent}. Explicitly, the result is the following
theorem.

\begin{thm}
  \label{cullershalen}
  Let $G=H_1*_MH_2$ be a graph of free groups such that each $H_i$
  has rank $2$. If $G\onto\free_3$ then $M$ is cyclic or trivial.
\end{thm}

\noindent One can derive upper bounds on the rank of the intersection given
lower bounds on the rank of the join. This has also been observed in
the nice article of Kent \cite{kent}, where some upper bounds are
explicitly computed. The proof of Theorem~\ref{cullershalen} presented
here differs only slightly from his. In the broadest terms, the two
articles share with most papers in the subject an analysis of
immersions of graphs, a method that dates back to Stallings
\cite{stallings0}. Specifically, Kent uses directly the topological
pushout of a pair of graphs along the core of their pullback, a graph
which appears here as the underlying graph of a reduced graph of
graphs.

\subsection*{Acknowledgements} The authors are very grateful to
Richard Kent for many discussions on this topic, in particular those
regarding Theorem~\ref{cullershalen}. The first author thanks the
California Institute of Technology for its hospitality during a visit
when this work began. The second author would like to thank Ben Klaff
for bringing the question of Culler and Shalen to his
attention. Finally, many thanks to the referee for several useful
comments and suggestions, especially a much simplified proof of
Lemma~\ref{lem:intersectlemma}.

%---------------------------------------------------------------------------
%---------------------------------------------------------------------------
\section{Graphs of graphs}
\label{sec:graphsofgraphs}

\noindent A \emph{graph of graphs} is a finite graph of spaces such that
all vertex spaces are combinatorial graphs and all edge maps are
embeddings. Below are some simple operations on graphs of graphs. All
vertices and edges are indicated by lower case letters, and their
associated spaces will be denoted by the corresponding letter in upper
case.  We will not keep track of orientation here despite its
occasional importance---we trust the reader to sort out this simple
matter when it arises.  Let $X$ be a graph of graphs with vertices
$v_i$ and edges $e_j$.

\begin{enumerate}
\item[({\bf M1})] \emph{Making vertex and edge spaces connected}: Let
  $V_{i,1},\dotsc,V_{i,n_i}$ be the connected components of the vertex
  space $V_i$ associated to the vertex $v_i$ of the underlying graph
  $G$, and $E_{j,1},\dotsc,E_{j,m_i}$ be the connected components of
  the edge space $E_j$. We construct a new graph of graphs as
  follows. First, we build the underlying graph. For each $i$ and $j$,
  we take a collection of vertices $v_{i,k}$ and edges $e_{j,l}$, one
  for each connected component of each vertex space and edge space,
  respectively. We label $v_{i,k}$ with $V_{i,k}$ and $e_{j,l}$ with
  $E_{j,l}$, and attach $e_{j,l}$ to $v_{i,k}$ if the image of
  $E_{j,l}$ in $V_i$ is contained in $V_{i,k}$. The attaching maps for
  this graph of graphs are the obvious ones. If $e_{j,l}$ is adjacent
  to $v_{i,k}$, then we attach an end of $E_{j,l}\times I$ to
  $V_{i,k}$ by the inclusion map.
\item[({\bf M2})] \emph{Removing unnecessary vertices}: If $V$ is a
  vertex space with exactly two incident edges such that both
  inclusions are isomorphisms, we remove $v$ and regard the pair of incident
  edges as a single edge. If $V$ has one incident edge and the
  inclusion is an isomorphism, we remove $V$ and the incident edge.
\item[({\bf M3})] \emph{Removing isolated edges}: If a vertex space
  $V$ has an edge $e$ that is not the image of an edge from an
  incident edge space, we remove $e$ from $V$.
\item[({\bf M4})] \emph{Collapsing free edges or vertices}: We call an
  edge $e$ of a vertex space $V$ \emph{free} if it is the image of
  only one edge from the collection of incident edge spaces, say
  $e'\subset E$.  In this case, we remove $e$ and $e'$ from $V$ and
  $E$. If a vertex space $V$ is a point and has only one incident edge
  space, we remove $v$ and the incident edge.
\end{enumerate}

\noindent A graph of graphs is \emph{reduced} if any application of
these operations leaves the space unchanged. Notice that any graph of
graphs can be converted to a reduced graph of graphs by greedily
applying \textrm{(M1)} through \textrm{(M4)}.  The remaining requisite
operations on graphs of graphs are blow ups and blow downs at a
vertex.

\begin{enumerate}
\item[({\bf M5})] For a vertex space $V$, divide the incident edge
  spaces into two classes $E_1,\dotsc,E_n$ and $E_{n+1},\dotsc,E_m$,
  and let $V_1$ ($V_2$, resp.) be the union of the images of $E_i$,
  $i\leq n$ ($i>n$, resp.). When $V_1 \cap V_2$ is non-trivial, we
  replace $V$ by $V_1\sqcup V_2$ and introduce a new edge $v_1\cap
  v_2$ with the edge graph $V_1\cap V_2$. Next, we attach $E_i$ to
  $V_1$ for $i\leq n$, $E_i$ to $V_2$ for $i>n$, and the newly
  introduced edge space $V_1\cap V_2$ to $V_1$ and $V_2$ via the
  inclusion maps.
\item[({\bf M6})] \emph{Blow up}: We blow up a vertex by applying
  \textrm{(M5)}. We pass to connected components of the newly created
  vertex and edge spaces via \textrm{(M1)}. Finally, we pass to the
  associated reduced graph of graphs using \textrm{(M2)}.
\item[({\bf M7})] \emph{Blow down}: Let $E$ be an edge space of a
  graph of graphs. If the two embeddings of $E$ have disjoint images,
  that is, $\iota(E)\cap\tau(E)=\emptyset$, then we remove the edge
  $e$ of the underlying graph and identify the two endpoints of
  $e$. Finally, the graph carried by the new vertex is the one
  obtained by identifying the vertex space(s) at the ends of $e$ by
  setting $\iota(f)=\tau(f)$, where $f$ is either a vertex or an edge
  of $E$.
\end{enumerate}

\begin{rem}
  Notice that if $X$ has no free or isolated edges, then the space
  obtained by blowing up a vertex with an application of \textrm{(M6)}
  also has no free or isolated edges. Also, when $V$ is connected, it
  follows that $V_1 \cap V_2$ is non-trivial and thus \textrm{(M5)} is
  applicable.
\end{rem}

\noindent The \emph{horizontal subgraph} of a graph of graphs is the
graph obtained by restricting vertex and edge spaces to vertices. The
\emph{mid-graph} of a graph of graphs is the graph obtained by
restricting vertex and edge spaces to midpoints of edges. These two
subgraphs are denoted $\Gamma_H$ and $\Gamma_M$, respectively. Note
that neither of these graphs is necessarily connected. If $X$ is
reduced, then $\Gamma_M$ and $\Gamma_H$ do not have any valence one
vertices. Conversely, if either one of them has a valence one vertex,
then there must be a free edge or vertex in $X$. If there are isolated
edges, then a component of $\Gamma_M$ is a point. In particular, if
$X$ is reduced, then every component of $\Gamma_M$ has nontrivial
fundamental group.

\begin{lemma}
  \label{lem:homotopytype}
  Blowing up and blowing down are homotopy equivalences.
\end{lemma}

\noindent This follows easily upon observing that if two of the edges
introduced during a blowup are both adjacent to a vertex introduced
during the blowup, then the images of the edge spaces they carry are
disjoint. That the remaining moves, other than \textrm{(M3)}, preserve
the homotopy types of $X$, $\Gamma_H$ and $\Gamma_M$, is clear. Note
that \textrm{(M3)} only serves to remove trivial components of
$\Gamma_M$.\smallskip\smallskip

\noindent It is important to know when to blow up $X$. The following
lemma achieves this.

\begin{lemma}
  \label{lem:intersectlemma}
  Let $\Delta$ be a connected graph with a collection
  $\mathcal{C}=\{\Delta_i\}_{i=1}^m$ of (not necessarily
  distinct) connected subgraphs. If every edge of $\Delta$ is
  contained in at least two $\Delta_i$ and $m>3$, then after
  relabeling the $\Delta_j$, there is a partition
  $\mathcal{C}_1=\set{\Delta_1,\dotsc,\Delta_n}$ and
  $\mathcal{C}_2=\set{\Delta_{n+1},\dotsc,\Delta_m}$ of $\mathcal{C}$
  such that at least two $\Delta_i$ in $\mathcal{C}_1$ intersect
  nontrivially and at least two $\Delta_i$ in $\mathcal{C}_2$
  intersect nontrivially.
\end{lemma}

\begin{proof}
  It suffices to find distinct $A,B,C,D\in\mathcal{C}$ such that
  $A\cap B\neq\emptyset$ and $C\cap D\neq \emptyset$. If all triple
  intersections are empty then $\mathcal{C}$ has at most two elements
  by connectivity of $\Delta$. Let $A,B,C\in\mathcal{C}$ such that
  $A\cap B\cap C\neq\emptyset$. Since $\Delta$ is connected, there is
  some $D$ meeting, again without loss, $C$.
\end{proof}

\begin{rem}
Notice that if $V$ is a vertex space of a
reduced graph of graphs $X$ with at least four incident edge spaces,
then we can use Lemma~\ref{lem:intersectlemma} to ensure that
\textrm{(M5)} is applicable.
\end{rem}

\noindent Let $X$ be a reduced graph of graphs such that all vertex and
edge spaces are connected. The space $X$ has an underlying graph that
we shall denote by $\Gamma_U(X)$. Let $m(X)$ be the highest valence of
vertex of $\Gamma_U(X)$, $n(x)$ the number of vertices with valence
$m(X)$, and $\chi(\Gamma_U(X))$ the Euler characteristic of
$\Gamma_U(X)$. The \emph{complexity} of $X$ is the lexicographically
ordered $3$--tuple
\[ c(X):=(\chi(\Gamma_U(X)),m(X),n(X)). \]
We call a blowup of a vertex $v$ using two sets of edge spaces
satisfying Lemma~\ref{lem:intersectlemma} \emph{nontrivial}. Our next
lemma justifies this terminology.

\begin{lemma}\label{complexity}
  Let $X$ be reduced and $m(X)>3$. If $X'$ is obtained from $X$ via a
  nontrivial application of \textrm{(M6)} to a vertex $v$ with valence
  $m(X)$, then $c(X')<c(X)$.
\end{lemma}

\begin{proof}
  Let $\set{v_i}$ be the vertices of $X'$ introduced during a blow up
  of $X$ at the vertex $v$. These vertices must have valence at least
  two, as otherwise $X$ has a free edge and is not reduced. We assume
  contrary to the claim that $c(X)=c(X')$. If the Euler
  characteristics of the underlying graphs of $X$ and $X'$ are equal,
  then the subgraph $B$ spanned by the edges associated to the
  connected components of $V_1\cap V_2$ must be a tree. First observe
  that it is connected as otherwise $V$ could not have been
  connected. Second, if $B$ is not a tree, then the Euler
  characteristic of the underlying graph must decrease. As $B$ is a
  tree we have
  \begin{equation}
    \label{1}
    1-\frac{1}{2}\textrm{valence}(v)=\sum_i \pr{1-\frac{1}{2}\textrm{valence}(v_i)}.
  \end{equation}
  If both $m(X')=m(X)$ and $n(X')=n(X)$, then all but one of the
  vertices $v_{i_0}$ has valence two since there are no valence one
  vertices making a positive contribution to the sum on the right hand
  side of (\ref{1}). Therefore, every component of $V_1$ (the
  alternative is handled identically) is the image of exactly one
  incident edge space from one element of the partition of edges
  incident to $v$. However, this is impossible since the blowup $X'$
  was assumed to be nontrivial.
\end{proof}

\noindent We are now ready to state our main result.

\begin{thm}\label{reducetovalencethree}
  Every graph of graphs $X$ such that no connected component of
  $\Gamma_M$ is a tree can be converted to a reduced graph of graphs
  $X'$ all of whose vertex groups have valence three. There is a
  homotopy equivalence
  $(X',\Gamma_H',\Gamma_M')\to(X,\Gamma_H,\Gamma_M)$.
\end{thm}

\noindent The \emph{corank} of a group $G$ is the maximal rank of a free
group that it maps onto and will be denoted by $\CR(G)$. Before
proving Theorem~\ref{reducetovalencethree}, a few remarks are in
order.  First, observe that if $X$ is reduced, then the natural map
$\pi_1(X)\to\pi_1(\Gamma_U(X))$ is surjective. Second, the complexity
of all graphs of graphs homotopy equivalent to $X$ is bounded below by
$(1-\CR(\pi_1(X)),3,0)$. That said, we now give a proof of
Theorem~\ref{reducetovalencethree}.

\begin{proof}[Proof of Theorem~\ref{reducetovalencethree}]
  First we apply $\textrm{(M4)}$ until there are no free edges. This
  does not change the homotopy type of the triple
  $(X,\Gamma_H,\Gamma_M)$. There are no isolated edges since each
  component of $\Gamma_M$ is assumed to have nontrivial fundamental
  group. Next, we pass to connected components of edge and vertex
  spaces and then pass to the associated reduced graph of graphs by
  removing valence two vertex spaces. Let $X$ be a reduced graph of
  graphs, and consider a sequence $\set{X_i}$ starting with $X$ such
  that $X_i$ is obtained from $X_{i-1}$ by nontrivially blowing up a
  maximal valence vertex. Since all the $X_i$ are homotopy equivalent
  and the maps $\pi_1(X_i)\to\pi_1(\Gamma_U(X_i))$, $i>0$, are
  surjective, $c(X_i)\geq (1-\CR(\pi_1(X)),3,0)$. According to
  Lemma~\ref{complexity}, $c(X_i)>c(X_{i+1})$. Since the complexity is
  bounded below, for some $n$, $X_n$ has only valence three vertices.
\end{proof}

\noindent A graph of graphs represents a graph of free groups when the
$\epsilon$--neighborhood of $\Gamma_M$ is a product
$I\times\Gamma_M$. In this case there are two natural immersions
$\Gamma_M\to\Gamma_H$ in the sense of Stallings
\cite{stallings0}. Moreover, there is an immersion
$\Gamma_H\to\Gamma_U$. We say such a graph of graphs is
\emph{representing}. Conversely, suppose that
$G=\Delta(H_1,\dotsc,H_k,M_1,\dotsb,M_l)$ is a graph of free groups
with vertex groups $H_i$, edge groups $M_j$, and that there is a map
$\gamma\co G\to\free$ which embeds each $H_i$. Let $\iota_j$ and
$\tau_j$ be the two inclusion maps $M_j\into H_{\iota(j)}$ and
$M_j\into H_{\tau(j)}$. Represent $\free$ as the fundamental group of
a marked labeled graph $R$ with one vertex, and find immersions of
marked labeled graphs $\eta_i\co\Gamma_{H_i}\to R$ representing
$\gamma\vert_{H_i}$, and $\mu_j\co\Gamma_{M_j}\to R$ representing
$\gamma\vert_{M_j}$. We choose the notation $\Gamma_{H_i}$ in
anticipation of the fact that they are the connected components of the
horizontal subgraph of the graph of graphs under construction.\smallskip\smallskip

\noindent The immersion $\mu_j$ factors through $\eta_{\iota(j)}$ and
${\eta_{\tau(j)}}$ via $\iota_j\co\Gamma_{M_j} \to
\Gamma_{H_{\iota(j)}}$ and $\tau_j\co\Gamma_{M_j} \to
\Gamma_{H_{\tau(j)}}$. We construct a space $X$ by taking the
$\Gamma_{M_j}\times I$ as edge spaces, taking the $\Gamma_{H_i}$ as
vertex spaces, and using as attaching maps
$\iota_j\co\Gamma_{M_j}\times\{0\}\to \Gamma_{H_{\iota(j)}}$ and
$\tau_j\co\Gamma_{M_j}\times\{1\}\to \Gamma_{H_{\tau(j)}}$. Let
$\alpha_j\co \Gamma_{M_j}\times I\to\Gamma_{M_j}$ be the projection
to the first factor. Since $\eta_{\iota(j)}\circ\iota_j=\mu_j$ and
$\eta_{\tau(j)}\circ\iota_j=\mu_j$ there is a well defined map
$\pi\co X\to R$ which restricts to $\eta_i$ and agrees with
$\mu_j\circ\alpha_j$.\smallskip\smallskip

\noindent We now endow $X$ with the structure of a graph of
graphs. Let $b$ be the base point of $R$. Let $V=\pi^{-1}(b)$ and
$E_l=\pi^{-1}(m_l)$, where $m_l$ is the midpoint of an edge $e_l$ of
$R$. Each edge $e_l$ of $R$ induces two maps of $E_l$ to $V$, each of
which is an embedding. That these are embeddings can be seen as
follows. If one fails to be injective on vertices of $E_l$, then some
$\Gamma_{H_i}\to R$ is not an immersion. If it is injective on
vertices but not on edges, then some $\Gamma_{M_j}\to R$ is not
immersed. Thus, we may use this data to endow $X$ with the structure
of a graph of graphs. By Theorem~\ref{reducetovalencethree}, we can
repeatedly blow up $X$ until we produce a graph of graphs $X'$ all of
whose vertices have valence three. If \textrm{(M3)} is ever applied in
the process, then it must be that some $M_i$ was
trivial.

\begin{rem}
Let $w$ be a vertex of a vertex space $V$
of a graph of graphs $X$. It follows that $w$ is a vertex of
$\Gamma_H$ and the valence of $w$ in $\Gamma_H$ is exactly the number
of edge graphs incident to $V$ whose images contain $w$. If $X$ is
reduced and $V$ has valence three, then there must be a vertex of $V$
which is contained in the image of all three incident edge
graphs. Moreover, if $\Gamma_M$ has a valence three vertex $w$, then
the images of $w$ in $\Gamma_H$ must each have valence three in
$\Gamma_H$.
\end{rem}

\begin{proof}[Proof of Theorem~\ref{cullershalen}]
  We begin by representing $G\onto\free_3$ by a map from a graph of
  graphs $X$ to a bouquet of three circles $R$. Note that since
  $G\onto \free_3$, at least one of the maps $H_i \to \free_3$ is injective.
  If the other fails to be injective, the result is immediate. Consequently, we are
  reduced to the case when both are injective and thus
  the construction above can be implemented. By blowing $X$ up, we may
  assume that $X$ has only valence three vertices. The map
  $G\onto\free_3$ factors through the map
  $G\cong\pi_1(X)\onto\pi_1(\Gamma_U)$, and the rank of $\Gamma_U$
  must be either $3$ or $4$. If the latter holds, then $M$ is trivial
  and the theorem holds. If $\Gamma_U$ has rank $3$, since all
  vertices of $\Gamma_U$ are valence three, there must be exactly
  four. By the remark above, if $\Gamma_M$ has a valence three vertex,
  then the map from the set of valence three vertices of $\Gamma_H$ to
  the set of valence three vertices of $\Gamma_U$ cannot be
  injective. Since the two components of $\Gamma_H$ each have
  fundamental group $\free_2$, they have $2$ valence three vertices
  apiece. However, this implies the contradiction that $\Gamma_U$ has
  at most $3$ valence three vertices. Thus, $\Gamma_M$ has vertices of
  valence at most two and so has rank at most one, as claimed.
\end{proof}

\begin{rem}
By \cite{louder}, $M$ is contained in the
subgroup generated by a basis element in at least one of $H_1$ or
$H_2$.
\end{rem}

\noindent Other inequalities of this type are easily obtained through
an analysis of a reduced valence three graph of graphs representing
the intersection. In particular, special cases of the Hanna Neumann
conjecture can be verified with this analysis. For explicit
inequalities, we refer the reader to Kent \cite{kent} who has also
derived them.

%---------------------------------------------------------------------------
%---------------------------------------------------------------------------
\section{Intersections of subgroups of free groups}

\noindent Let $H_1$ and $H_2$ be subgroups of a fixed free group
$\free$. If $G=\Delta(H_1,H_2;M_j)$, a graph of free groups with two
vertex groups $\{H_i\}$, edge groups $\{M_j\}$, with no monogons and a
a map $\pi\co G\onto\free$ embedding each of the factors $H_i$, then
the vertex spaces of a graph of graphs $X$ representing $\Delta$ are
bipartite.\smallskip\smallskip

\noindent A graph of graphs is \emph{simple-edged} if no vertex space has a
bigon. To relate reduced graphs of graphs to intersections of free
groups we need to understand what happens when a graph of graphs $X$
as above is not simple-edged. Let $p$ and $q$ be the midpoints of a
pair of offending edges, $\Gamma=\Gamma_U(X)$ the underlying graph of
$X$, and give the edges of $\Gamma$ distinct oriented labels. The
labeling of $\Gamma$ induces labelings of $\Gamma_M$ and
$\Gamma_{H_i}$. Let $\Gamma_M'$ be the labeled graph obtained by
identifying $p$ and $q$. By folding the labeled graph $\Gamma_M'$ (see
for instance \cite{stallings0}), we obtain a labeled graph $\Gamma_K$
with fundamental group $K=\pi_1(\Gamma_K,p)$. Folding endows
$\Gamma_K$ with a pair of immersions
$\nu_i\co\Gamma_K\to\Gamma_{H_i}$. In addition, there is an immersion
$\eta_j\co\Gamma_{M_j}\to\Gamma_K$ and the edge map
$\Gamma_{M_j}\to\Gamma_{H_i}$ is just $\nu_i\circ\eta_j$.\smallskip\smallskip

\noindent We must consider two cases with regard to $p$ and $q$ after
folding. Namely, the midpoints $p$ and $q$ are either in the same
component of $\Gamma_M$ or in distinct components of $\Gamma_M$.  We
address the latter first and assume, without loss of generality, that
$\Gamma_{M_1}$ and $\Gamma_{M_2}$ are the components of $\Gamma_M$
containing $p$ and $q$.  Compute the fundamental groups of
$\Gamma_{H_1}$ and $\Gamma_{H_2}$ with respect to the images of $p$
(which coincide with the images of $q$). From this we see that
$\pi(H_1)\cap\pi(H_2)$ contains $\pi(M_1)$ and $\pi(M_2)$. If
$\pi(M_1)\not<\pi(M_2)$ and $\pi(M_2)\not<\pi(M_1)$, then each
inclusion $\pi(M_i)\into \pi(H_1)\cap\pi(H_2)$ is proper and the image
of the fundamental group of $\Gamma_K$, computed with respect to the
image of $p$, is precisely $\langle M_1,M_2\rangle$. If neither
$\eta_1$ nor $\eta_2$ is an isomorphism of labeled graphs, then $K$
properly contains $M_1$ and $M_2$. In the event we are in the first
case, without loss of generality, we shall assume that $p,q$ are
contained in $\Gamma_{M_1}$. We identify the vertices $p$ and $q$ of
$\Gamma_M$ and then fold to obtain a labeled graph $\Gamma_K$. As
before, the immersion $\Gamma_{M_1}\to\Gamma_{H_i}$ factors through
the induced immersion $\Gamma_K\to\Gamma_{H_i}$. In this case
$\Gamma_M\to\Gamma_K$ cannot be an isomorphism of graphs and $H_1\cap
H_2$ properly contains $M_1$.\smallskip\smallskip

\noindent Let $X$ be a reduced simple-edged graph of graphs with
underlying graph $\Gamma=\Gamma_U(X)$. Let $\mathcal{X}(\Gamma)$ be
the collection of reduced simple-edged graphs of graphs with
underlying graph $\Gamma$. If $X,X'\in\mathcal{X}(\Gamma)$, then
$X\leq X'$ if there is a map of graphs of spaces $X\to X'$ such that
all restrictions to vertex and edge spaces are embeddings. We can
restrict to the subcollection $\mathcal{X}_{X}(\Gamma)$ such that for
each $X'\in\mathcal{X}_{X}(\Gamma)$ there is a map $X\to X'$ and the
map $\Gamma_H(X)\to\Gamma_H(X')$ is a graph isomorphism. Clearly
$\mathcal{X}_{X}(\Gamma)$ has a maximal element $Y$. To link reduced
simple-edged graphs of graphs to the strengthened Hanna Neumann
conjecture, we only need to observe that since $X$ is simple-edged,
each component of $\Gamma_M(X)$ is an embedded subgraph of
$\Gamma_M(Y)$ (i.e.\ the fundamental groups of components of
$\Gamma_M(X)$ are free factors of the respective components of
$\Gamma_M(Y)$).\smallskip\smallskip

\noindent The strengthened Hanna Neumann conjecture then implies that
if $G$ is as above and the associated graph of graphs is simple-edged,
then
\[\chi(\Gamma_{H_1})\chi(\Gamma_{H_2})+\chi(\Gamma_{M})\geq 0\]
The equivalence of the amalgamated graph conjecture and the
strengthened Hanna Neumann conjecture of \cite{dicks} follows
immediately from the observation that the vertex and edge spaces of a
representing simple-edged graph of graphs can be written as in the
statement of Dicks' theorem. We leave the details of the construction
of this correspondence to the reader, though we state a version of the
equivalence for completeness.\smallskip\smallskip

\noindent Let $X$ be a simple-edged reduced graph of graphs all of whose
vertices are valence three that represents a homomorphism
$\Delta(H_1,H_2,M_j)\to\free$. Let $v_i$ be the vertices of
$\Gamma_U(X)$, and for each $i$, let $\Delta_i$ be the intersection of
the images of the three edge spaces incident to $v_i$. Finally, let $\Delta$ be
the disjoint union of the $\Delta_i$,
$\Sigma_1=\Delta\cap\Gamma_{H_1}$ and
$\Sigma_2=\Delta\cap\Gamma_{H_2}$, and $\mu$ be the number of edges in
$\Delta$.

\begin{thm}
  \[\chi(H_1)\chi(H_2)+\sum_i\chi(M_i)=\frac{1}{4}\vert\Sigma_1\vert\dot\vert\Sigma_2\vert-\frac{1}{2}\mu\]
\end{thm}

\noindent The proof is straightforward. The number of valence three
vertices of $\Gamma_{H_i}$ is $\vert\Sigma_i\vert$, $\Gamma_{H_i}$ has
only valence two or valence three vertices, and the Euler
characteristic of $H_i$ is therefore
$-\frac{1}{2}\vert\Sigma_i\vert$. The Euler characteristic of each
$\Gamma_{M_j}$ is computed in the same manner. $\mu$ is the number of
valence three vertices of $\Gamma_{M}$. In this formulation, the
amalgamated graph conjecture simply states that if one is given a
reduced simple-edged representing graph of graphs whose horizontal
graph has two components, then the right hand side of the above
equality is nonnegative.

%---------------------------------------------------------------------------
%---------------------------------------------------------------------------

%---------------------------------------------------------------------------
%---------------------------------------------------------------------------

\end{document}